\documentclass[a4paper,11pt]{article}

\usepackage[english]{babel}

\usepackage{amssymb}
\usepackage{latexsym}
\usepackage{amsfonts}
\usepackage{amsmath} 
\usepackage{epsfig}
\usepackage{graphics}

\newcommand{\T}{{\mathbb T}}

\newcommand{\R}{{\mathbb R}}
\newcommand{\C}{{\mathbb C}}

\newcommand{\Z}{{\mathbb Z}}

\newcommand{\gl}{\operatorname{{\mathfrak gl}}}
\newcommand{\GL}{\operatorname{GL}}
\newcommand{\Ad}{\operatorname{Ad}}
\newcommand{\ad}{\operatorname{ad}}
\newcommand{\im}{\operatorname{im}}
\newcommand{\codim}{\operatorname{codim}}
\newcommand{\orb}{\cO}

\newcommand{\KAM}{\textsc{kam}\ }
\newcommand{\lcu}{\textsc{lcu}}
\newcommand{\BHTa}{\textsc{bht}{\rm (i)}}
\newcommand{\BHTb}{\textsc{bht}{\rm (ii)}}

\newcommand{\cB}{{\mathcal{B}}}

\newcommand{\cO}{{\mathcal{O}}}
\newcommand{\cX}{{\mathcal{X}}}
\newcommand{\id}{\operatorname{Id}}
\newcommand{\diag}{\operatorname{diag}}
\newcommand{\I}{{\rm i}}
\newcommand{\sg}{\Sigma}

\numberwithin{equation}{section}

\newenvironment{remarks}{\vskip 8pt \noindent {\bf Remarks.}\rm \begin{itemize}}{\end{itemize}}

\newtheorem{theorem}{Theorem}
\newtheorem{definition}[theorem]{Definition}

\newtheorem{corollary}[theorem]{Corollary}
\newtheorem{lemma}[theorem]{Lemma}

\newtheorem{Ex}{Example}
\newenvironment{ex}{\begin{Ex}\rm}{\hspace*{\fill}$\Box$\end{Ex}}

\newenvironment{mylist}{%
\begin{list}{$\bullet$}{\setlength{\itemsep}{3pt}
\setlength{\parsep}{0pt} \setlength{\topsep}{3pt}
\setlength{\partopsep}{0pt}
\setlength{\labelwidth}{0.8cm}
\setlength{\leftmargin}{1.0cm}
\setlength{\labelsep}{0.2cm}}}{\end{list}}

\date{23 September 2008}

\title{Quasi-periodic stability of normally resonant tori.}
  
\author{Henk W. Broer\thanks{Dept.\ of Mathematics and Computing Science,
        University of Groningen, POBox 407,
        9700~AK~Groningen, The Netherlands} ,  
        M.~Cristina Ciocci\thanks{Dept.\ of Mathematics,
        Imperial College London, Campus South Kensington, Queen's Gate,
        London, SW72AZ, UK} , 
        Heinz Han{\ss}mann\thanks{Mathematisch Instituut,
        Universiteit Utrecht, Postbus 80.010,
        3508~TA~Utrecht, The Netherlands} \mbox{} and
        Andr\'e Vanderbauwhede\thanks{Dept.\ of Pure Mathematics and
        Computeralgebra, University of Gent,
        Krijgslaan~281, 9000~Gent, Belgium} }

\begin{document}

\maketitle

\begin{abstract}
\noindent
We study quasi-periodic tori under a normal-internal resonance, possibly
with multiple eigenvalues.
Two non-degeneracy conditions play a role.
The first of these generalizes invertibility of the Floquet matrix and
prevents drift of the lower dimensional torus.
The second condition involves a Kolmogorov-like variation of the internal
frequencies and simultaneously versality of the Floquet matrix unfolding.
We focus on the reversible setting, but our
results carry over to the Hamiltonian and dissipative contexts.
\end{abstract}

\bigskip
{\textbf{MSC-class:}} 37J40

\section{Introduction}

Persistence results for quasi-periodic motions were first proved for
maximal tori in Hamiltonian systems and became known as
Kolmogorov--Arnol'd--Moser(KAM)theory.
In~\cite{Mos 67} this was extended to lower dimensional tori and to
other contexts like volume preserving and reversible systems.
The r\^ole of the `modifying terms' in terms of system parameters
was clarified in~\cite{BHT, H} and the R\"ussmann
condition~\cite{BHS, rue01} allows to subsequently reduce the high
number of parameters to the bare minimum.

\medskip\noindent
These results yield what is called quasi-periodic (or normal linear)
stability, i.e.\ families of invariant tori persist under sufficiently
small perturbations when restricted to certain (measure-theoretically
large) Cantor sets.
The theorems in~\cite{BHS} make the crucial assumption that all
eigenvalues of the matrix~$\Omega$ describing the normal linear
behaviour be simple.
This implies in particular that $\det \Omega \neq 0$ (except for the
disspative case and the high-dimensional volume preserving case, where
this condition is explicitly added).
Multiple resonances are admitted in~\cite{BHN, CMC2, Hoo} and the aim
of the present paper is to admit zero eigenvalues without weakening
the conclusion of quasi-periodic stability.

\subsection{Setting and results}
\label{settingandresults}

We work on the phase space
$M=\T^n \times \R^m \times \R^{2p}$, where $\T^n=(\R/2\pi\Z)^n$ is the $n$-torus
on which we use coordinates $x=(x_1, \ldots, x_n) \pmod{2\pi }$,
while on $\R^m$ and $\R^{2p}$ we use respectively $y=(y_1,\ldots,y_m)$ and
$z=(z_1,\ldots,z_{2p})$. In such coordinates a vector field on $M$ takes the form
\[
   \dot{x}= f(x,y,z), \qquad 
   \dot{y}=g(x,y,z), \qquad
   \dot{z}=h(x,y,z),
\]
or in vector field notation:
\begin{eqnarray}
   X(x,y,z) =f(x,y,z)\partial_x + g(x,y,z)\partial_y + h(x,y,z)\partial_z. 
\label{e:vctXfgh}
\end{eqnarray}
We assume that the vector field $X$ depends analytically on all variables,
including possible parameters which we suppress for the moment; referring
to~\cite{BHT,H,Po82} we note that our results remain valid when
`analyticity' is replaced by `a sufficiently high degree of differentiability'.
An invariant torus $T$ of a vector field $X$ is called parallel if a smooth
conjugation exists of the restriction 
$X|_T$ with a constant vector field $\dot{x} = \omega$ on $\T^n$. 
The vector $\omega = (\omega_1,\omega_2,\ldots,x_n)\in \R^n$
is the (internal) frequency vector of $T$.
The parallel torus is quasi-periodic when the frequencies are independent over the rationals.

\medskip\noindent
We are concerned with persistence of quasi-periodic tori under small
perturbations, and to fix thoughts we concentrate\footnote{We give explicit formulations 
for reversible vector fields, but the results remain valid for e.g.\ dissipative, Hamiltonian or
volume-preserving systems (vector fields and maps), where equivariance is also optional.} 
on the reversible setting. To define reversibility we consider an involution (i.e.\ $G^2=I$)
\begin{equation}
\label{e:G}
   G: {M} \longrightarrow {M}, \quad (x,y,z)\mapsto (-x,y,Rz),
\end{equation}
with $R\in \GL(2p, \R)$ a linear involution on $\R^{2p}$ such that 
\[
   \dim \textrm{Fix}(R) = \dim \left\{z\in\R^{2p} \mid Rz=z \right\} = p.
\]
The vector field $X$ is then called $G$-reversible
(or reversible for short) if 
\[
   G_*(X)=-X.
\] 
Using (\ref{e:vctXfgh}) this reversibility condition takes the explicit form
\begin{eqnarray*}
   f(-x,y,Rz) & = &    \phantom{-} f(x,y,z),  \nonumber \\
   g(-x,y,Rz) & = &  -  g(x,y,z), \\
   h(-x,y,Rz) & = &  - Rh(x,y,z), \nonumber 
\end{eqnarray*}
valid for all $(x,y,z)\in M$.

\medskip\noindent
Following \cite{BH,BHS,BHT,H} the vector field $X$ is called integrable if it is equivariant with
respect to the group action
\begin{displaymath}
   \T^n \times M \longrightarrow M, \quad
   (\xi, (x, y, z)) \mapsto (\xi + x, y, z)
\end{displaymath}
of $\T^n$ on $M$, or in other words, if the functions $f$, $g$ and $h$
in~\eqref{e:vctXfgh} are independent of the $x$-variable(s).
Such an integrable vector field 
\begin{equation}
\label{vf:Xintegrable}
   X(x,y,z)=f(y,z)\partial_x + g(y,z)\partial_y + h(y,z)\partial_z
\end{equation} 
is reversible if
\begin{equation}
\label{integr+revers}
   f(y,Rz)=f(y,z), \qquad g(y,Rz)=-g(y,z) \qquad \textrm{and} \qquad
   h(y,Rz)= - Rh(y,z)
\end{equation}
for all $(y,z)\in \R^m \times \R^{2p}$; this implies $g(y,z)=0$ for all
$(y,z)\in \R^m \times \textrm{Fix}(R)$. 
In case $h(0, 0) = 0$ the\footnote{Often one has a whole family 
$T_y = \T^n \times \{ y \} \times \{ 0 \}$ of invariant tori. While we 
are especially interested in bifurcations, the variable $y$ will still 
act as a parameter, now unfolding the bifurcation scenario.} 
$n$-torus $T_0 = \T^n \times \{ 0 \} \times \{ 0 \}$
is invariant under the flow of the vector field $X$.
The normal linear part $N(X)$ of~\eqref{integr+revers} at~$T_0$ is given by 
\begin{equation}
\label{N-X}
   N(X)(x,y,z) = \omega\partial_x + \Omega\, z\,\partial_z, 
\end{equation}
with
\[
   \omega= f(0,0) \quad\textrm{and} \quad\Omega = D_z h(0,0).
\]
We denote the subspace of infinitesimally reversible linear operators on $\R^{2p}$ by 
$\gl_-(2p;\R)$ and by $\gl_+(2p;\R)$ the subspace of all $R$-equivariant
linear operators on $\R^{2p}$, i.e.\ 
\[
   \gl_\pm(2p;\R) = \{ \Omega \in \gl(2p;\R) \mid \Omega R = \pm R \Omega \}.
\]
In order to define the non-degeneracy of~\eqref{vf:Xintegrable}
at the invariant torus $T_0$ we consider the subspaces
\begin{displaymath}
   \cX_{lin}^{\pm G} = \left\{ \omega \partial_x + \Omega z \partial_z \mid
   \omega \in \R^n, \Omega \in \gl_{\pm}(2p;\R) \right\}
\end{displaymath}
of the spaces $\cX^{-G}$ of all $G$-reversible vector fields on~$M$ and
$\cX^{+G}$ of all $G$-equivariant vector fields, satisfying $G_*(X) = +X$.
For $X \in \cX^{-G}$ the adjoint operator
\begin{displaymath}
   \ad N(X) : \cX \longrightarrow \cX,
   \quad Y \mapsto [ N(X), Y ]
\end{displaymath}
maps $\cX^{\pm G}$ into $\cX^{\mp G}$; a similar statement is true for~$\cX_{lin}^{\pm G}$.

\medskip\noindent
Our interest concerns purely $G$-reversible vector fields,
and $G$-reversible vector fields that are furthermore equivariant with
respect to
\begin{equation}
\label{elldeck}
   F_l : M \longrightarrow M, \quad
   (x, y, z) \mapsto
   (x_1 - \frac{2 \pi}{l}, x_*,
   y, z_I, {\rm e}^{ \frac{2 \pi \I}{l}}\, z_{II}).
\end{equation}
Here $z_{II} \cong z_{2j-1} + \I z_{2j}$ singles out two of the
$z$-variables in a complexified form and
$z_I = ( z_1, z_2, \ldots, z_{2j-2}, z_{2j+1},\ldots, z_{2p} )$
contains the remaining $z$-variables.
To allow for a unified formulation of our results we define a reversing symmetry group $\Sigma$
and a character (a group homomorphism) $\chi:\Sigma \longrightarrow \{\pm 1\}$ as follows:
\begin{mylist}
\item[(i)] In the purely reversible case we set $\Sigma:=\{\id,G\}$ and $\chi(G):=-1$.
\item[(ii)] In the equivariant-reversible case we define $\Sigma$ as the group generated by
$G$ and $F_l$ and define $\chi$ by $\chi(G):=-1$ and $\chi(F_l):=1$.
\end{mylist}
In both cases $\Sigma$ is isomorphic to $\Z_2\ltimes Z_l$, the dihedral group of order~$2l$.
When $l = 1$ the generator $F_1 = \id$ of course is superfluous. 
For both cases we put
\begin{eqnarray*}
   \cX^+ & = & \left\{ X \in \cX \mid
   \mbox{$E_*(X) = X$ for all $E \in \sg$} \right\}  \\
   \cX^- & = & \left\{ X \in \cX \mid
   \mbox{$E_*(X) = \chi(E) X$ for all $E \in \sg$} \right\}
\end{eqnarray*}
together with $\cX_{lin}^{\pm} = \cX_{lin}^{\pm G} \cap \cX^{\pm}$.
Furthermore we let $\cB^+$ and~$\cB^-$ consist of the constant vector fields
in $\cX^+$ and~$\cX^-$, respectively and denote by
\[
   \orb(\Omega_0) = \left\{ \Ad(A)\cdot \Omega_0 := A \,\Omega_0 \, A^{-1}
   \mid A \in \GL_+(2p;\R) \right\}
\]
the orbit under the adjoint action of $\GL_+(2p;\R)$ on $\gl_-(2p;\R)$.

\begin{definition}[Broer, Huitema and Takens \cite{BHT}]
\label{definition1}
The parametrised\,\footnote{The r\^ole of the external parameter~$\lambda$ occurring 
in Definition~\ref{definition1} can be (partially) taken by the internal parameter~$y$.} 
vector field $X_{\lambda}$ with linearization
$N(X_{\lambda})(x, y, z) =
 \omega(\lambda) \partial_x + \Omega(\lambda) z \partial_z$
is \emph{non-degenerate} at $\lambda=\lambda_0\in\R^s$ if
\begin{mylist}
\item[\BHTa\phantom{i}] $\ker \ad N({X}_{\lambda_0}) \cap \cB^+ = \{ 0 \}$;
\item[\BHTb]  at $\lambda=\lambda_0$ the mapping
$(\omega, \Omega) : \R^s \longrightarrow \R^n \times \gl_-(2p;\R),
 \lambda \mapsto (\omega(\lambda),\Omega(\lambda))$
is transverse to
$\{\omega(\lambda_0)\} \times \orb\left(\Omega(\lambda_0)\right)$.
\end{mylist}
\end{definition}

\noindent
The two non-degeneracy conditions \BHTa\ and~\BHTb\ generalize the
condition that $\ad N({X}_{\lambda_0})$ has to be invertible, a
requirement that lies at the basis of Mel'nikov's conditions
(\eqref{eq:diophantine} with $|\ell| \neq 0$).
One also speaks of BHT~non-degeneracy.
Compared to the formulation in~\cite{BHT}, \S~8a2 the requirement
that $\Omega(\lambda_0)$ have only simple eigenvalues is dropped.
The extension to multiple normal frequencies was developed
in~\cite{BHN,CMC2,Hoo} for invertible~$\Omega(\lambda_0)$; we
return to the original formulation of~\BHTa.

\medskip\noindent
To formulate the strong non-resonance condition necessary for persistence
of invariant tori we introduce for $\Omega\in \gl_-(2p;\R)$ the normal
frequency mapping $\alpha : \gl_-(2p;\R) \longrightarrow \R^{2p}$ 
where the components of $\alpha(\Omega)$ 
are equal to the imaginary parts of the eigenvalues of
$\Omega\in \gl_-(2p;\R)$. 
Higher multiplicities are taken into account 
by repeating each eigenvalue as many times as necessary.  

\begin{definition}
\label{diophantine-cond}
A pair $(\omega, \Omega)\in \R^n \times \gl_-(2p;\R)$ is said to satisfy
a \emph{Diophantine condition} if there exist constants $\tau > n-1$ 
and $\gamma >0$ such that
\begin{equation}
\label{eq:diophantine}
   \mid \langle k, \omega \rangle + \langle \ell, \alpha(\Omega) \rangle \mid
   \geq \gamma |k|^{-\tau}
\end{equation}
for all $k \in \Z^n \setminus \{0\}$ and
$\ell \in \Z^{2p}$ with $|\ell| \leq 2$. 
\end{definition}

\noindent
This condition is independent of the way in which we have ordered the
components of $\alpha(\Omega)$; also, if $(\omega,\Omega)$
satisfies~\eqref{eq:diophantine} then the same is true for all
$(\omega, \widetilde{\Omega})$ with $\widetilde{\Omega} \in \orb(\Omega)$.
For each $\Gamma \subset P$ we define the associated Diophantine subset
\begin{displaymath}
   \Gamma_{\gamma}:= \left\{ \lambda\in \Gamma \mid
   (\omega(\lambda), \Omega(\lambda)) \,
   \mbox{satisfies~\eqref{eq:diophantine}} \, \right\}.
\end{displaymath}
When $\Gamma$ is a small neighbourhood of some $\lambda_0\in P$
where $X$ is non-degenerate
then $\Gamma_{\gamma}$ is nowhere dense but with large measure 
(provided that $\gamma$ is sufficiently small).

\begin{theorem}
\label{main-thm}
Let $X \in \cX^-$ be a family of $\sg$-reversible integrable vector fields
that is non-degenerate at $\lambda_0 \in P$. 
Then there exists $\gamma_0 > 0$ such that for all $0 < \gamma < \gamma_0$
the following is true. There exists a neighbourhood $\Gamma$ of $\lambda_0$,
neighbourhoods $\mathcal{Y}$ and $\mathcal{Z}$ of the origin in 
respectively $\R^m$ and $\R^{2p}$,
and a neighbourhood $\mathcal{U}$ of $X$ in the compact-open topology 
on~$\cX^-$ such that for each $Z\in \mathcal{U}$ one can find a mapping 
$\Phi:\T^n \times \mathcal{Y} \times \mathcal{Z} \times \Gamma
 \longrightarrow M\times P$ 
of the form
\[
   \Phi(x,y,z,\omega,\mu) =
   \left( x + \widetilde{U}(x,\omega,\mu), y+ \widetilde{V}(x,y,\omega,\mu),
   z + \widetilde{W}(x,y,z,\omega,\mu), \omega + \widetilde{\Lambda}_1(\omega,\mu),
   \mu + \widetilde{\Lambda}_2(\omega,\mu) \right)
\]
for which the following holds.
\begin{mylist}
\item[\textup{(i)}] 
	The mapping $\Phi$ is $\sg$-equivariant, real-analytic in the $x$-variable 
        and normally affine in the $y$ and $z$ variables.
\item[\textup{(ii)}]
	The mapping $\Phi$ is $C^{\infty}$-close to the identity  and is a $C^{\infty}$-diffeomorphism 
        onto its image.
\item[\textup{(iii)}]
	The restriction of $\Phi$ to the Cantor set
	$\T^n \times \{ 0 \}\times \{ 0 \}  \times \Gamma_{\gamma}$
	of Diophantine $X$-invariant tori conjugates $X$ to $Z$. 
	The restiction of $\Phi$ to 
	$\T^n \times \mathcal{Y} \times \mathcal{Z} \times \Gamma_{\gamma}$
	also preserves the normal linear behaviour to these invariant tori.
\end{mylist}
\end{theorem}

\noindent
In terms of~\cite{BHS,BHT}, the conclusion of Theorem \ref{main-thm} expresses
that the family $X$ is quasi-periodically stable,
i.e., structurally stable on a union of (Diophantine) quasi-periodic tori.
This allows to condense Theorem~\ref{main-thm} to the statement
that non-degenerate $\sg$-reversible integrable vector fields
are quasi-periodically\footnote{In \cite{BHN} one speaks of
`normal linear stability' instead.} stable.
Quasi-periodic stability implies that for every small perturbation~$Z$
there exists a $Z$-invariant `Cantor set' $V\subset M\times P$
which is a $C^{\infty}$-near-identity diffeomorphic image of the foliation 
$\T^n \times \{ 0 \}\times \{ 0 \}\times \Gamma_{\gamma}$ of $n$-tori. 
In the tori this diffeomorphism is an analytic conjugacy from $X$ to~$Z$,
which also preserves the normal linear behaviour.

\subsection{Normal-internal resonances}
\label{normalinternalresonances}

Resonances are at the core of the problems one has to solve when trying to
prove quasi-periodic stability -- persistence of elliptic invariant tori
\begin{displaymath}
    T_y = \T^n \times \{y\} \times 0 \subseteq N:= \T^n \times \R^m \times \R^{2p}
\end{displaymath} 
under small perturbation. 
The strong non-resonance conditions \eqref{eq:diophantine} exclude in fact
four types of resonances.
An internal resonance
\begin{displaymath}
   \langle k, \omega \rangle = 0
   \quad \mbox{for some $0 \neq k \in \Z^n$}
\end{displaymath}
prevents the parallel flow on $T_y $
to have a dense orbit whence the invariant torus is not a (minimal)
dynamical object, but rather the union of closed invariant subtori.
One cannot expect such an $n$-torus to persist, cf.~\cite{Moser, sev86},
for the same reason that a circle consisting of equilibria breaks up
under perturbation (generically with only finitely many equilibria in the perturbed
system).
Such resonances are excluded by~\eqref{eq:diophantine} when taking $\ell = 0$.

\medskip\noindent
For $|\ell| = 1$ the inequalities~\eqref{eq:diophantine} constitute
the first Mel'nikov condition, cf.~\cite{bou97, mel65, XuYo01},
and concern the normal-internal resonances
\begin{equation}
\label{foldMelnikov}
   \langle k, \omega \rangle = \alpha_j \quad
   \mbox{with fixed $k \in \Z^n$ and $j \in \{ 1, \ldots, p \}$.}
\end{equation}
Passing to co-rotating co-ordinates on~$N$ yields this resonance
with $k=0$, cf.~\cite{BHJVW03, BrVe92}.
This is a $2$-step procedure.
First $k$ is brought into the form $k = (k_1, 0, \ldots, 0)$ by
means of a preliminary transformation
\begin{equation}
\label{fiveman}
   N \longrightarrow N , \quad (x, y, z) \mapsto (\sigma x, y, z)
\end{equation}
with $\sigma \in SL(n, \Z)$.
For the second step we again write
$z_I = (z_1, z_2, \ldots, z_{2j-2}, z_{2j+1},\ldots, z_{2p})$,
$z_{II} = (z_{2j-1},z_{2j})$ and complexify $z_{II} \cong z_{2j-1} + \I z_{2j}$.
Then we perform a Van der Pol transformation
\begin{equation}
\label{interpol}
   N \longrightarrow N , \quad (x, y, z) \mapsto (x, y, z_I, {\rm e}^{\I k_1 x_1} z_{II}).
\end{equation}
The transformed vector field has a vanishing
normal frequency $\alpha_j = 0$.
Hence, already constant perturbations
\begin{displaymath}
   \beta \partial_z = \beta_{2j-1} \partial_{z_{2j-1}}
   + \beta_{2j} \partial_{z_{2j}},
   \quad \beta_{2j-1}, \beta_{2j} \in \R
\end{displaymath}
make the tori~$T_y$ move in a way that cannot be compensated on the
linear level.
Condition \BHTa\ in Definition~\ref{definition1} prevents such perturbations
whence Theorem~\ref{main-thm} yields quasi-periodic stability, see also
Corollary~\ref{corrollary-3} in Section~\ref{perturbationproblem}.
An alternative to this condition is to take (generic) higher order terms of
the unperturbed vector field into account.
This typically results in bifurcation scenarios that turn out to be
quasi-periodically stable (in an appropriate sense) as well,
cf.~\cite{BBH, BHJVW03, h1h98, h1h07, wag05}.

\medskip\noindent
The remaining possibility $|\ell| = 2$ in~\eqref{eq:diophantine} excludes the
normal-internal resonances
\begin{equation}
\label{hopfMelnikov}
   \langle k, \omega \rangle = \alpha_i \pm \alpha_j \quad
   \mbox{with fixed $k \in \Z^n$ and $i \neq j \in \{ 1, \ldots, p \}$}
\end{equation}
and
\begin{equation}
\label{flipMelnikov}
   \langle k, \omega \rangle = 2 \alpha_j \quad
   \mbox{with fixed $k \in \Z^n$ and $j \in \{ 1, \ldots, p \}$.}
\end{equation}
For~\eqref{hopfMelnikov} one can again achieve $k=0$ in co-rotating
co-ordinates, cf.~\cite{XuYo01}, turning this normal-internal
resonance into the normal resonance
\begin{displaymath}
   0 \neq \alpha_i = \pm \alpha_j, \quad
   i \neq j \in \{ 1, \ldots, p \} .
\end{displaymath}
While now the invertibility of~$\Omega$ does yield quasi-periodic
stability of~$T_y$, see~\cite{BHN,CMC2} and Corollary~\ref{corrollary-1}
in Section~\ref{perturbationproblem}, the normal behaviour still may be
drastically affected.
Using the $y$-variable as a parameter, e.g.\ the normal linear matrix
\begin{displaymath}
   \left(
   \begin{array}{cccc}
      0 & 1 & 0 & 0  \\
      -1 & 0 & 0 & 0  \\
      0 & 0 & 0 & -1  \\
      0 & 0 & 1 & 0
   \end{array}
   \right)
\end{displaymath}
in a conservative setting unfolds (or deforms) both to elliptic and to
hyperbolic behaviour.
Here it is the bifurcation scenario involving the surrounding tori of
dimension $n+1$ and $n+2$ that can only be captured by taking higher
order terms of the unperturbed vector field into account; quasi-periodic
stability was achieved in~\cite{BCH, BHH07, h1h07,Hoo} for the simplest
conservative bifurcation scenarios.

\medskip\noindent
The remaining case~\eqref{flipMelnikov} is meaningful only if not already
implied by~\eqref{foldMelnikov}, so assume that~\eqref{eq:diophantine} holds
with $|\ell| \leq 1$.
Then we can still achieve $k=0$ in co-rotating co-ordinates, but now on
a $2$-fold covering $M \longrightarrow N$ defined as follows.
The preliminary transformation~\eqref{fiveman} brings the
resonance vector $k$ into the form $k = (k_1, 0, \ldots, 0)$
with $k_1$ odd.
The Van der Pol transformation is no longer a mapping from~$N$
to itself, but a covering mapping from
$M = \T^n \times \R^m \times \R^{2p}$ onto~$N$ defined by
\[
   \Pi: M\longrightarrow N, \, \, \ (x_1,x_*,y,z)\mapsto 
   (2 x_1, x_*, y, z_{I},  {\rm e}^{\I k_1 x_1} z_{II}).
\]
Here $x_1 \in \T^1$ and $x_* = (x_2,\ldots,x_n)\in \T^{n-1}.$
The deck group $\Z_2 = \{\id,F\}$ of this $2$-fold covering is
generated by the involution
\begin{equation}
\label{deck}
   F: M\longrightarrow M, \, \, \ (x_1,x_*,y ,z_I, z_{II})\mapsto 
   (x_1 - \pi,x_*,y, z_I, -z_{II}).
\end{equation}
This means that
\[
   \Pi \circ F = \Pi.
\]
Note that this is the special case $l=2$ of~\eqref{elldeck}, the
corresponding quasi-periodic stability is stated in Corollary~\ref{corrollary-2}
of Section~\ref{perturbationproblem}.
The resulting frequency halving (or quasi-periodic period-doubling)
bifurcation scenarios are quasi-periodically stable in the
dissipative~\cite{BBH} and Hamiltonian~\cite{h1h07} settings and 
similarly reversible frequency-halving bifurcations may be expected
to occur if appropriate non-degeneracy conditions on the higher order
terms are fulfilled.

\medskip\noindent
In~\cite{bou97, XuYo01}\footnote{These papers consider Hamiltonian systems,
but we expect the results to carry over to the reversible context.} the
second Mel'nikov condition (\eqref{eq:diophantine} with $|\ell| = 2$) is avoided
completely, i.e.\ also simultaneous normal-internal resonances 
\eqref{hopfMelnikov} and~\eqref{flipMelnikov} with
differing $k \in \Z^n$ are allowed.
The price to pay for this approach is that any control on the linear
behaviour is completely lost.
For instance, double eigenvalues $\pm \I \alpha_1 = \pm \I \alpha_2$
generically unfold to a Krein collision, where an elliptic torus
evolves a $4$-dimensional normal direction of focus-focus type.
Such changes cannot be captured without persistence of the (normal)
linear behaviour.

\subsection{Contents and conclusions}
\label{contentsandconclusions}

This paper fits in the framework of parametrised
KAM~theory~\cite{BBH, BCH, BHH07, BHJVW03, BHS, h1h07}
that originates from Moser \cite{Mos 67};
in fact we present a generalization of~\cite{BHN, BH, BHT},
as well as of \cite{CMC2,Hoo,H}.
In the next section we explicitly work out two examples to which
Theorem~\ref{main-thm} applies.
In Section~\ref{perturbationproblem} we elaborate the conditions
of Theorem~\ref{main-thm} and also formulate three corollaries;
the Corollaries~\ref{corrollary-2} and~\ref{corrollary-3} make
the novel results of this paper explicit.
The proof of Theorem~\ref{main-thm} is sketched in Section~\ref{section4}.
The necessary unfolding theory, which plays a key r\^ole, is deferred to the Appendix.

\medskip\noindent
Our approach allows for normal-internal resonances~\eqref{hopfMelnikov}
and~\eqref{flipMelnikov} (and possibly also~\eqref{foldMelnikov}) with $k \in \Z^n$ fixed.
The ensuing deformations of the linear behaviour coming from the perturbation
are taken care of by considering a versal unfolding of the linear part
$\dot{z} = \Omega z$ of the unperturbed vector field, i.e., an unfolding
that already contains all possible deformations.
The necessary parameters are provided by $y \in \R^m$; the possibility
that $m \geq n$ distinguishes the reversible context from the Hamiltonian setting.
An alternative is to let the system depend on external parameters~$\lambda$,
where variation of $(y, \lambda)$ versally unfolds the linear part.

\medskip\noindent
The proof in Section~\ref{section4} is formulated in terms of filtered Lie
algebras and therefore exceeds the reversible setting, carrying over to
other contexts that can be formulated in these terms, notably the dissipative,
volume preserving and Hamiltonian contexts; possibly combined with equivariance,
cf.~\cite{BHN,BHS}.
In the Hamiltonian case this answers a conjecture formulated in~\cite{h1h07}
to the positive.
For dissipative systems this has already been used in~\cite{BBH} when
proving quasi-periodic stability of the frequency-halving bifurcation
scenario.
We expect that appropriate higher order terms in~\eqref{unperturbed}
allow to obtain a similar result for reversible systems.

\medskip\noindent
The unfolding~\eqref{unf_0_geomult_1} recovers the result for the case $p=2$
that was obtained in \cite{Iooss}. 
There a $4$-dimensional reversible system with a 	
codimension $2$ singularity at the origin is studied 
by formal normal forms together with the persistence 
of the associated codimension $1$ bifurcation phenomena. 
It would be interesting to investigate the persistence of the corresponding 
bifurcation scenario in the {\sc kam} setting. 
Note that an additional $F$-equivariance next to the $G$-reversibility would
enforce the origin to be an equilibrium for the entire non-linear family, an
assumption that is made in~\cite{Iooss}.

\paragraph{Acknowledgment. } The second author received financial support by the European Community's 6th Framework Programme, Marie Curie Intraeuropean Fellowship EC contract Ref.\ MEIF-CT-2005-515291, award Nr.\ MATH P00286.\\

\section{Applications}
\label{applications}

We illustrate our results with two examples that explicitly show how our
assumptions enter and what extra conclusions can be drawn.

\begin{ex}\label{exampleone} \textbf{(Quasi-periodic response solutions)}
To show how to check the appropriate assumptions we consider the simple
example of a $1$-parameter family of quasi-periodically forced
oscillators
\begin{equation}
\label{responseproblem}
   \ddot{z} = f_{\mu}(t, z, \dot{z}) = h_{\mu}(t, \omega t, z, \dot{z}),
\end{equation}
with a fixed frequency $\omega,$ for instance we take
$\omega = \frac{1}{2}(\sqrt{5} - 1)$ (the golden mean number).
The forcing~$h_\mu$ is $2\pi$-periodic in the first two arguments.
The search is for quasi-periodic response solutions with this same
frequency vector $(1, \omega)$.

\medskip\noindent
Putting $z_1 = z$, $z_2 = \dot{z}$ we can rewrite~\eqref{responseproblem}
as an autonomous system
\begin{eqnarray*}
   \dot{x}_1 & = & 1  \\
   \dot{x}_2 & = & \omega  \\
   \dot{z}_1 & = & z_2  \\
   \dot{z}_2 & = & h_{\mu}(x, z) = \bar{h}_{\mu}(z) + \tilde{h}_{\mu}(x, z)
\end{eqnarray*}
on $\T^2 \times \R^2$ where we split~$h_{\mu}$ into the
average $\bar{h}_{\mu}$ over~$\T^2 \times \{ z \}$ and the oscillating
part $\tilde{h}_{\mu} = h_{\mu} - \bar{h}_{\mu}$.
The integrable vector field $X_{\mu}$ given by
\begin{eqnarray*}
   \dot{x}_1 & = & 1  \\
   \dot{x}_2 & = & \omega  \\
   \dot{z}_1 & = & z_2  \\
   \dot{z}_2 & = & \bar{h}_{\mu}(z)
\end{eqnarray*}
has invariant $2$-tori for all $z_1 \in \R$ with $\bar{h}_{\mu}(z_1, 0) = 0$.
These correspond to response solutions of the forced oscillator.

\medskip\noindent
Note that we allowed for $h_{\mu}$ to depend explicitly on~$z_2$ whence
$z_2 \mapsto -z_2$ is not a reversing symmetry.
We impose the system to be reversible with respect to
\begin{displaymath}
   (x_1, x_2, z_1, z_2) \mapsto (-x_1, -x_2, - z_1, z_2),
\end{displaymath}
in particular $\bar{h}_{\mu}$ depends on $z_1$ only through~$z_1^2$ and we
concentrate on the invariant torus at $z=0$.
The dominant part
\begin{displaymath}
   N(X_{\mu}) = \partial_{x_1} + \omega \partial_{x_2} + 
   \Omega(\mu) z \partial_z
\end{displaymath}
has the parameter-dependent $2 \times 2$ matrix
\begin{displaymath}
   \Omega(\mu) = \left(
   \begin{array}{cc}
      0 & 1 \\
      \partial_1 \bar{h}_{\mu}(0) &
      \partial_2 \bar{h}_{\mu}(0)
   \end{array} \right)
\end{displaymath}
which is invertible whenever $\partial_1 \bar{h}_{\mu}(0) \neq 0$.
However, the non-degeneracy condition \BHTa\ is also fulfilled
if $\partial_1 \bar{h}_{\mu}(0) = 0$ since the eigenvector to
the resulting eigenvalue~$0$ is not invariant under the involution
\begin{displaymath}
   R = \left(
   \begin{array}{cc}
      -1 & 0 \\  0 & 1
   \end{array} \right) .
\end{displaymath}
From this we conclude that condition \BHTa\ is always satisfied. 

\medskip\noindent 
The non-degeneracy condition \BHTb\ is satisfied when
\begin{equation}\label{check}
   \frac{{\rm d}}{{\rm d} \mu} \partial_1 \bar{h}_{\mu}(0) \neq 0 .
\end{equation}
Thus, the system is BHT non-degenerate as soon as \eqref{check} holds true. 
Therefore, given this by Theorem~\ref{main-thm} (for an explicit formulation
see Corollary~\ref{corrollary-3}), if the oscillating part $\tilde{h}$ is
sufficiently small, the forced oscillator~\eqref{responseproblem} has a response solution
near $z=0$, with linear behaviour changing where $\partial_1 \bar{h}_{\mu}(0)$
passes through zero. 
\end{ex}

\begin{remarks}
\item[(i)] Earlier for the existence of a response solution as an extra requirement
           the condition $ \partial_1 \bar{h}_{\mu}(0) \neq 0$
           was needed~\cite{BHN,BH,BHS,mel65,Mos 67}.
\item[(ii)] The stability change of the response solution as $\partial_1 \bar{h}_{\mu}(0)$
            passes through zero leads in the periodic case to additional periodic solutions
            bifurcating off from $z=0$, cf.~\cite{lam94, LaCa95, PCQW90}.
            We expect such bifurcations to carry over to the quasi-periodic case.
\end{remarks}

\noindent
We now return to the setting of the introduction where the normal-internal
resonance~\eqref{flipMelnikov} led to a perturbation problem on a
$2$:$1$~covering space.  
The next example shows how the normally linear vector fields
on the covering and the base-space relate to one another. 

\begin{ex}\label{exampletwo}\textbf{(Multiple normal-internal resonance)}
On the phase space
\[
N = \T^2 \times \R^3 \times \R^4 = \{x, y, z\}
\]
we consider the normally linear vector field
\begin{displaymath}
   Y = 2 \partial_{x_1} + \omega \partial_{x_2} + \Omega(\mu) z \partial_z
\end{displaymath}
with
\begin{displaymath}
   \Omega(\mu) = \left(
   \begin{array}{cccc}
      0 & -1 - \mu_1 & 1 & 0  \\
      1 + \mu_1 & 0 & 0 & 1  \\
      - \mu_2 & 0 &  0 & -1 - \mu_1  \\
      0 & - \mu_2 & 1 + \mu_1 & 0
   \end{array} \right)
\end{displaymath}
where we think of the parameters $\nu = (\omega, \mu) \in \R^3$ as been
obtained from $y \in \R^3$ by localization~\eqref{localisation}.
The eigenvalues $\pm \I (1 + \mu_1) \pm \sqrt{-\mu_2}$ of~$\Omega(\mu)$
yield at $\mu = 0$ the normal frequency $\alpha = \pm \I$ that has two
normal-internal resonances \eqref{hopfMelnikov} and~\eqref{flipMelnikov}
with the same $k = (1, 0) \in \Z^2$.
Complexifying both $\zeta_I \cong \zeta_1+\I \zeta_2$ and
$\zeta_{II} \cong \zeta_3 + \I \zeta_4$ on the covering space
\begin{displaymath}
   \hat{N} = \R/(4\pi\Z) \times \T \times \R^3 \times \R^4
    = \{ \xi_1, \xi_2, \eta, \zeta \}
\end{displaymath}
we have the covering mapping
\begin{displaymath}
   \Pi: \hat{N} \longrightarrow N, \, \, \ (\xi_1, \xi_2, \eta, \zeta) \mapsto 
   (\xi_1 \mbox{\rm mod} (2\pi \Z), \xi_2, \eta, \diag [{\rm e}^{\frac{1}{2}\I \xi_1}]\, \zeta).
\end{displaymath}
This leads to the deck transformation 
\begin{equation}
\label{doubledeck}
   F: \hat{N} \longrightarrow \hat{N}, \, \, \
   (\xi_1, \xi_2, \eta, \zeta) \mapsto (\xi_1 - 2\pi, \xi_2, \eta, -\zeta)
\end{equation}
and the lifted vector field
\begin{displaymath}
   \hat{Y} = \hat{\omega}_1 \partial_{\xi_1} + \omega_2 \partial_{\xi_2}
    + \hat{\Omega}(\mu) \zeta \partial_{\zeta}
\end{displaymath}
on~$\hat{N}$ satisfying $\Pi_* \hat{Y} = Y$.
In this setting $\dot{\xi}_1 = \dot{x}_1,$ implying that $\hat{\omega}_1 = 2$ 
and the corresponding periods are $\hat{T}_1 = 2\pi$ and $T_1 = \pi$,
so $\hat{T}_1 = 2 T_1$ as should be expected. 

\medskip\noindent
Regarding the Floquet matrices $\Omega$ and $\hat{\Omega}$ we have
\begin{eqnarray*}
\dot{z} & = & \diag[\frac{1}{2} \I \dot{\xi}_1 {\rm e}^{\frac{1}{2}\I \xi_1}]\,\zeta  + 
              \diag[{\rm e}^{\frac{1}{2}\I \xi_1}]\,\dot{\zeta} \\
        & = & \diag[{\rm e}^{\frac{1}{2}\I \xi_1}]\,
              \left( \frac{1}{2} \I \dot{\xi}_1 \zeta + \hat{\Omega} \zeta \right) \\
        & = & \diag[{\rm e}^{\frac{1}{2}\I \xi_1}]\,  \left( \I \id + \hat{\Omega}\right) \zeta \\
        & = & \diag[{\rm e}^{\frac{1}{2}\I \xi_1}]\,  \left( \I \id + \hat{\Omega}\right)  
              \diag[{\rm e}^{- \frac{1}{2}\I \xi_1}] \, z.
\end{eqnarray*}
Apparently
\[
\Omega = 
\diag[{\rm e}^{\frac{1}{2}\I \xi_1}] \left(\I \id + \hat{\Omega}\right)
\diag[{\rm e}^{-\frac{1}{2}\I \xi_1}]
= \I \id + \hat{\Omega},
\]
and the resulting family
\begin{displaymath}
   \hat{\Omega}(\mu) = \Omega(\mu) - \I \id = \left(
   \begin{array}{cccc}
      0 & -\mu_1 & 1 & 0  \\
      \mu_1 & 0 & 0 & 1  \\
      - \mu_2 & 0 &  0 & -\mu_1  \\
      0 & - \mu_2 & \mu_1 & 0
   \end{array} \right)
\end{displaymath}
of matrices is the \lcu\ of  
\begin{displaymath}
   \hat{\Omega}(0) = \left(
   \begin{array}{cccc}
      0 & 0 & 1 & 0  \\
      0 & 0 & 0 & 1  \\
      0 & 0 & 0 & 0  \\
      0 & 0 & 0 & 0
   \end{array} \right).
\end{displaymath}
Every perturbation of~$Y$ on~$N$ can be lifted to a perturbation of~$\hat{Y}$
on~$\hat{N}$ that respects the deck transformation~\eqref{doubledeck} and rescaling time 
we can always arrange $\dot{x}_1 = 2,$ i.e.\ that the first frequency equals $2$.  
Applying Theorem~\ref{main-thm} (for an explicit formulation see Corollary~\ref{corrollary-2})
we may conclude that~$\hat{Y}$ is quasi-periodically stable and this implies quasi-periodic
stability of~$Y$.
\end{ex}

\noindent
It should be noted that such an application of \KAM Theory goes beyond
the possibilities of  \cite{BHN,CMC2,Hoo}. 
For $n = 1$ the full bifurcation scenario has been addressed in~\cite{BrFu93} and it 
would be interesting to develop the extension for periodic to quasi-periodic orbits.

\section{Main results}
\label{perturbationproblem}

In the perturbation problem we work on the phase space 
$M = \T^n\times\R^m\times\R^{2p} = \{x,y,z\}$ where we are dealing with 
a `dominant part'
\begin{equation}
\label{dp}
   \dot{x} = \omega, \quad
   \dot{y} = 0, \quad
   \dot{z} = \Omega z, \quad
   \mbox{ or } \quad
   X = \omega\partial_x + \Omega z \partial_z 
\end{equation}
in vector field notation.
While it is always possible to translate a single given invariant torus to
$T_0 = \T^n \times \{ 0 \} \times \{ 0 \}$, it is an assumption on the system
that this torus can be embedded in a whole family
$T_y = \T^n \times \{ y \} \times \{ 0 \}$ of invariant tori
parametrised by~$y$.
This can be equivalently stated as
\begin{equation}
\label{familyoftori}
   h(y, 0) = 0 \quad \mbox{for all $y \in \R^m$,}
\end{equation}
and the non-degeneracy condition \BHTa\ in Definition~\ref{definition1}
ensures that this assumption can be justified.
For each $\epsilon>0$  the scaling operator
\begin{equation}
\label{scaling-op}
   \mathcal{D}_\epsilon : M \longrightarrow M, \; (x,y,z) \mapsto 
   \left(x, \frac{y}{\epsilon} , \frac{z}{\epsilon^2}\right)
\end{equation}
commutes with $G$ and with the $\T^n$-action on $M$, and hence preserves
reversibility and integrability.
Using~\eqref{vf:Xintegrable} and the linearity of $\mathcal{D}_\epsilon$
the push-forward  $\left(\mathcal{D}_\epsilon\right)_\ast\! X $ of $X$
under $\mathcal{D}_\epsilon$ takes the form
\begin{eqnarray*}
   \left(\mathcal{D}_\epsilon\right)_\ast\! X (x,y,z) & =  &
   \mathcal{D}_\epsilon \left( X\!\! \left(\mathcal{D}_\epsilon^{-1}(x,y,z)
   \right)\right) \\
   & = & f( \epsilon y, \epsilon^2 z) \partial_x
   + \frac{1}{\epsilon}g( \epsilon y, \epsilon^2 z) \partial_y
   + \frac{1}{\epsilon^2}h( \epsilon y, \epsilon^2 z) \partial_z.
\end{eqnarray*}
We can use~\eqref{integr+revers} to find that
$N(X):= {\displaystyle \lim_{ \epsilon\rightarrow 0 }}
 \left(\mathcal{D}_\epsilon\right)_\ast\! X$
is the dominant part~\eqref{N-X} of~$X$.
The vector field $N(X) = \omega \partial_x + \Omega z \partial_z$ 
is again reversible and integrable;
it is characterised by the frequency vector
$\omega=( \omega_1, \ldots, \omega_n)\in \R^n$ which describes the flow along the 
invariant tori $T_y$,
and by the matrix $\Omega\in \gl(2p;\R)$ which determines the linear flow
in the $z$-direction normal to the family of invariant tori.

\medskip\noindent
Since $\Omega$~does not depend on the angular variable $x \in \T^n$
the vector field $N(X)$ is in normal linear Floquet form.
The Floquet matrix $\Omega$ is {infinitesimally reversible}, satisfying
$\Omega R = - R \Omega$ because of the reversibility of the vector field~$X$.
Observe that if $\mu\in\C$ is an eigenvalue of $\Omega\in \gl_-(2p;\R)$
then so is $-\mu$.
Hence, the eigenvalues of $\Omega$ can be grouped into
complex quartets, conjugate purely imaginary pairs $\pm \I \alpha$,
symmetric real pairs and the eigenvalue zero with even algebraic multiplicity.

\subsection{Non-degeneracy conditions}
\label{Sec:NDs}

Since $\GL_+(2p;\R)$ is algebraic it follows that the orbit
$\orb(\Omega_0)$ is a smooth submanifold of $\gl_-(2p;\R)$.
The tangent space at $\Omega_0$ to this orbit is given by 
\begin{equation}
\label{tangent-to-orbit}
   T_{\Omega_0}\orb(\Omega_0) =
   \{\ad(A)\cdot \Omega_0 = A \Omega_0 - \Omega_0 A \mid A \in \gl_+(2p;\R) \}
   = \ad(\Omega_0)\left( \gl_+(2p;\R) \right),
\end{equation}
where we have used the fact that $\ad(A)\cdot \Omega =- \ad(\Omega) \cdot A$ for all
$A, \Omega \in \gl(2p;\R)$. An {unfolding} of $\Omega_0$ is a smooth mapping
\[
   \Omega: \R^s \longrightarrow \gl_-(2p;\R), \; \mu \mapsto \Omega(\mu)
\]
such that $\Omega(0)=\Omega_0$.
An unfolding is {versal} if it is transverse to $\orb(\Omega_0)$
at $\mu=0$, which requires that $s \geq \codim \orb(\Omega_0)$;
a versal unfolding with the minimal number of parameters (i.e.\ with $s$ equal
to the codimension of $\orb(\Omega_0)$ in $\gl_-(2p;\R)$) is called {miniversal}.
Using the Implicit Function Theorem it is easily seen that given a miniversal
unfolding $\Omega: \R^s \longrightarrow \gl_-(2p;\R)$ of $\Omega_0\in \gl_-(2p;\R)$
we can write each $\widetilde{\Omega}\in \gl_-(2p;\R)$ near
$\Omega_0$ in the form $\widetilde{\Omega} = \Ad(A)\cdot \Omega(\mu)$ for some
$(A,\mu)\in \gl_+(2p;\R) \times \R^s$ close to $(\id, 0)$
and depending smoothly on $\widetilde{\Omega}$.
In case all eigenvalues of~$\Omega_0$ are different from each other a
miniversal unfolding amounts to simultaneously deforming the eigenvalues,
see~\cite{BHT}.
Our approach yields persistence results independent of the eigenvalue
structure of~$\Omega_0$ (see \cite{Ba} for some other step towards
such general persistence results).
For more details on versal, miniversal (or universal) unfoldings we refer to \cite{A71,A83,Gib}.

\medskip\noindent
Property \BHTa\ generalizes the invertibility condition required
in the definition of non-degeneracy as it was formulated
in~\cite{BHN,BH, CMC2,Hoo}.
What is really needed for the proofs is the invertibility of the
linear operator
\begin{equation}
\label{revadjoint}
   \ad N({X}_{\lambda_0}) : \cB^+ \longrightarrow \cB^- 
\end{equation}
and since $\dim \textrm{Fix}(R)=\dim \textrm{Fix}(-R)$
this is fully captured by~\BHTa.
Computing
\begin{equation}
\label{betaadjoint}
   \ad N({X}_{\lambda_0}) (\beta \partial_z) =
   -{\Omega}(\lambda_0) \beta \partial_z
\end{equation}
shows that this certainly holds true if
$\det \Omega(\lambda_0)\neq 0$.
However, the condition \BHTa\ can still be satisfied if
$\det \Omega(\lambda_0)=0$, for example when
$\ker(\Omega(\lambda_0))\subset \textrm{Fix}(-R)$.
The Floquet matrix $\Omega(\lambda_0)$ may have zero
eigenvalues as long as the corresponding eigenvectors
do not lie in~$\cB^+$.

\begin{remarks}
\item[(i)] Up to now, the condition $\det \Omega_0 \neq 0$ was one of the
central assumptions for normal linear stability in the general dissipative context as well as 
in the volume preserving, symplectic and reversible contexts. 
Replacing this condition by \BHTa\ allows to extend the known theorems 
to the singular case of eigenvalue zero. 
\item[(ii)] Property \BHTa\ is persistent under small 
          variation of $\lambda$ near $\lambda_0$ because of 
          the upper-semi-continuity of the mapping 
          $\lambda \mapsto \rm{dim} \ker\Omega(\lambda).$
\end{remarks}

\noindent
Property \BHTb\ means that locally the frequency vector $\omega(\lambda)$
varies diffeomorphically with $\lambda,$ while `simultaneously' the local family
$\lambda \mapsto \Omega(\lambda)$ is a versal unfolding of $\Omega(\lambda_0)$ 
in the sense of \cite{A71,A83}. For earlier usage of this method 
in reversible \textsc{kam} Theory, see~\cite{BHN, BH, CMC2}.
In the Appendix we develop an appropriate
versal unfolding that depends linearly on $\lambda.$

\medskip\noindent
When trying to answer the persistence problem for~$T_y$ it is convenient
to focus on (a sufficiently small neighbourhood of) each of the invariant tori
$T_\nu$ ($\nu\in\R^m$) separately, considering the label $\nu\in \R^m$ of the
chosen torus as a parameter; 
formally this can be done by a localizing transformation, setting
\begin{equation}
\label{localisation}
   y= \nu + y_{loc}
   \qquad \textrm{and} \qquad
   X_{loc}(x, y_{loc}, z; \nu) := X(x, \nu + y_{loc}, z).
\end{equation}
This way we get a parametrised family of reversible and integrable
vector fields, still on the same state space $M$; in this localized
situation we concentrate on the persistence in a small neighbourhood of the
invariant torus $T_0$, corresponding to $(y_{loc},z)=(0,0)$.
For simplicity we absorb the additional parameter $\nu$ with the other
parameters which we may have and we also drop the subscript `$loc$'.

\medskip\noindent
The non-degeneracy condition~\BHTb\ requires that
$s \geq n + \codim \orb\left(\Omega(\lambda_0)\right)$;
in case all parameters originate from a localization procedure
this means that we should have
$m \geq n + \codim \orb\left(\Omega(\lambda_0)\right)$.
Assume now that $X_{\lambda}$ is non-degenerate at $\lambda_0\in\R^s$, and let
$(\omega_0,\Omega_0):= (\omega(\lambda_0), \Omega(\lambda_0))$.
Using a re-parametrisation and a parameter-dependent linear transformation in
the $z$-space we may assume that the parameter $\lambda$ takes the form
$\lambda=(\omega, \mu, \tilde{\mu})$ and belongs to a neighbourhood $P$ of
$\lambda_0:=(\omega_0,0,0)$ in $ \R^n \times \R^c \times \R^{s-n-c}$, while
the dominant part of the vector field reads
\begin{equation}
\label{unpert-vf-standard}
   N(X)(x,y,z,\omega,\mu,\tilde{\mu}) =
   \omega \partial_x + \Omega(\mu) z \partial_z,
\end{equation}
where $\Omega: \R^c \longrightarrow \gl_-(2p;\R)$ is a given miniversal unfolding
of~$\Omega_0$.
The $\tilde{\mu}$-part of the parameter does not appear in this expression for
the (unperturbed) vector field $X$; although it might appear explicitly in the
perturbations it turns out that $\tilde{\mu}$ plays no role at all in the
further analysis.
Therefore we suppress $\tilde{\mu}$ and just keep the essential
parameters $(\omega,\mu)$ and set $s=n+c$, with
$c=\codim \orb(\Omega(0))$.
The question of a particular choice for the miniversal unfolding
$\Omega(\mu)$ appearing in~\eqref{unpert-vf-standard}
is addressed in the Appendix.

\subsection{Consequences}
\label{Sec:mr}

We are given a family of integrable vector fields
\begin{equation}
\label{unperturbed}
   X(x, y, z, \omega, \mu) = 
   \left[ \omega + f(y, z, \omega, \mu) \right] \partial_x +
   g(y, z, \omega, \mu) \,\partial_y
   + \left[ \Omega(\mu) z + h(y, z, \omega, \mu) \right] \partial_z
\end{equation}
on the product $M \times P$ of phase space
$M = \T^n \times \R^m \times \R^{2p}$ and parameter
space $P \subseteq \R^s = \R^n \times \R^c$ with reversing symmetry
group $\sg$ generated by \eqref{e:G} and~\eqref{elldeck}.
For $l=1$ the latter is just the identity, but for $l \geq 2$ the
composition
\begin{equation}
\label{e:Hell}
   H_l := F_l \circ G : M \longrightarrow M, \quad
   (x, y, z) \mapsto (\frac{2 \pi}{l}-x_1, x_*, y, S_l R z),
\end{equation}
is another reversing symmetry and one may also characterise the
vector fields in~$\cX^-$ as being reversible with respect to the
two mappings $G$ and~$H_l$.
Note that in this characterisation $H_l$ may be replaced by
$F_l^i \circ G$ for any $i$ relative prime to~$l$.

\medskip\noindent
The coefficient functions $f, g$ and~$h$ entering~$X$ are higher order terms
in~$z$, satisfying
$f(y, 0, \omega, \mu) = g(y, 0, \omega, \mu) = h(y, 0, \omega, \mu)
 = D_z h(y, 0, \omega, \mu) = 0$
for all $y \in \R^m$, $(\omega, \mu) \in P$.
Within $\cX^-$ we consider perturbations $Z$ of~$X$ and write
\begin{displaymath}
   Z(x,y,z,\omega,\mu) = 
   \left[\omega + \tilde{f}(x,y,z,\omega,\mu)\right]\partial_x +
   \tilde{g}(x,y,z,\omega,\mu)\,\partial_y +
   \left [\Omega(\mu) z + \tilde{h}(x,y,z,\omega,\mu)\right] \partial_z  ;
\end{displaymath}
here the coefficient functions $\tilde{f}, \tilde{g}$ and~$\tilde{h}$ may
contain lower order terms but are close to $f, g$ and~$h$, respectively.
In this situation Theorem~\ref{main-thm} allows to conjugate $Z$ to~$X$
as far as Diophantine tori are concerned.

\medskip\noindent
The condition that $\Phi$ be a full conjugation from $X$ to $Z$ means that
$\Phi_*(X) = Z$, or equivalently $\left(\Phi^{-1}\right)_*(Z) = X$.
What is actually proved is the existence of a local diffeomorphism $\Phi$ such that
\begin{equation}
\label{actual-conjugation-cond}
   \left(\Phi^{-1}\right)_*(Z)(x,y,x,\omega,\mu) = N(X)(x,y,z,\omega,\mu)
   + O(|y|,|z|)\partial_x + O(|y|, |z|^2) \partial_y + O(|y|, |z|^2)\partial_z
\end{equation}
for all $(\omega,\mu)\in P_{\gamma}$ which are sufficiently close to $(\omega_0,0)$.
The property (\ref{actual-conjugation-cond}) implies that for all parameter values 
$(\omega,\mu)$ in the indicated Cantor set the $X$-invariant torus 
$\T^n \times \{0\} \times \{0\}$ is mapped by $\Phi$ into a $Z$-invariant torus 
on which the $Z$-flow is conjugate to the constant flow $\omega \partial_x$ on $\T^n$. 
This means that a Cantor subset of large measure of the family 
$\T^n \times \{0\}\times \{ 0 \} \times P$ of $X$-invariant tori 
survives the perturbation to $Z$.
The preservation of the normal linear behaviour means that the
normal linear vector fields $N(X)$ and~$N(Z)$ along two corresponding
invariant tori are conjugated by the derivative of the
$C^{\infty}$-near-identity diffeomorphism.

\medskip\noindent
In comparison to earlier results on persistence of lower-dimensional tori
the condition that all eigenvalues be simple is dropped in
Theorem~\ref{main-thm} and the condition $\det \Omega(0) \neq 0$ is
weakened to~\BHTa. Indeed, we have the following corollary.

\begin{corollary}[Ciocci \cite{CMC2}, Broer, Hoo and Naudot \cite{BHN}]
\label{corrollary-1}
Let the family $X \in \cX^-$ of $G$-reversible integrable vector fields
satisfy the non-degeneracy condition~\BHTb\ at
$\lambda_0 = (\omega_0, 0) \in P$, with $\Omega(0)$ invertible.
Then $X$ is quasi-periodically stable.
\end{corollary}

\noindent
Next to the above purely reversible case $l=1$ also the case $l=2$ of
a reversing symmetry group $\sg = \{ \id, F, G, H \}$ merits an
explicit formulation.
Here $H = H_2$ is given by~\eqref{e:Hell} and  
yields
\begin{eqnarray}\label{T-reversibility}
        f(y,S Rz)=f(y,z),
        \qquad g(y,S Rz)=-g(y,z)
        \qquad \mbox{and} \qquad h(y,S Rz)=-SR h(y,z)
\end{eqnarray}
for integrable vector fields where $S(z_I, z_{II}) = (z_I, -z_{II})$.
From (\ref{equivariance}) it follows that $f,g$ are even in $z_{II}$, 
while $h$ is odd in $z_{II}$. 
Moreover, (\ref{integr+revers}) and (\ref{T-reversibility}) imply that 
$g(y,z)=0$ for all $(y,z)\in \R^m\times \textup{Fix}(R)$ 
and also for all $(y,z)\in \R^m\times \textup{Fix}(S R)$. 

\begin{corollary} \label{corrollary-2}
Let $X \in \cX^-$ be a family of $G$-reversible $F$-equivariant
integrable vector fields that satisfies the non-degeneracy condition~\BHTb\ at
$\lambda_0 = (\omega_0, 0) \in P$, with $\Omega(0)$ invertible on~$\textup{Fix}(S)$.
Then $X$ is quasi-periodically stable.
\end{corollary}

\noindent
Again we allow for multiple eigenvalues, in particular the eigenvalue~$0$
may have multiplicity larger than two.
A similar statement holds in case of equivariance with respect
to~\eqref{elldeck} instead of~\eqref{deck}.

\medskip\noindent
In the covering setting of Section 1, we observe that the lift of an
integrable vector field again is integrable. 
In fact, if $\hat{X}$ is the lift to ${M}$ of an integrable vector field 
${X}$ on $N$, then $\Pi_*(\hat{X})={X}$ and $F_*\hat{X}=\hat{X}$. 
In case the second Mel'nikov condition is violated by a
resonance~\eqref{flipMelnikov} we can apply Corollary~\ref{corrollary-2}
on a $2$:$1$~covering space.
In Example \ref{exampletwo} of Section \ref{applications} we do this for
a double normal-internal resonance with fixed resonance vector $k \in \Z^2$.

\begin{corollary} \label{corrollary-3}
Let $X \in \cX^-$ be a family of $G$-reversible integrable vector fields that
satisfies the non-degeneracy condition~\BHTb\ at
$\lambda_0 = (\omega_0, 0) \in P$.
If $\ker \Omega(0)$ is contained
in $\textup{Fix}(-R)$ then $X$ is quasi-periodically stable.
\end{corollary}

\noindent
If $\ker \Omega(0) \subseteq \textup{Fix}(+R)$ we generically expect a
quasi-periodic centre-saddle bifurcation to take place, cf.~\cite{h1h98}.
Here violations of the first Mel'nikov condition prevents persistence of the
corresponding tori if not embedded in an appropriate bifurcation scenario.
The scaling \eqref{scaling-op} also can be applied to non-integrable systems,
making the non-integrable higher order terms a small perturbation.
It is then not automatic that the resulting dominant part is in
Floquet form, this is a necessary extra requirement that can be
thought of as generalization of integrability under which
quasi-periodic stability can still be achieved.
For a more thorough discussion of these questions see~\cite{BHT}.

\section{Sketch of proof}
\label{section4}

The proof of Theorem~\ref{main-thm} follows~\cite{BB, BHN, BHT} almost
verbatim (see also \cite{BH, CMC2, Hoo}) and here we just concentrate
on the novel aspects. 
The quite universal set-up of~\cite{BHT,Mos 67} is based on a Lie algebra approach,
using a standard Newtonian linearization procedure.
The conjugation $\Phi$ between the integrable and the perturbed family
is produced as the limit of an infinite iteration process.
The central ingredient of the proof is the solution of the linearized
problem, the so-called homological equation. 
The structure at hand, that is, the reversible symmetry group~$\sg$, is
phrased in terms of the Lie algebras $\cX^{\pm}$, $\cX_{lin}^{\pm}$
and~$\cB^{\pm}$ and therefore automatically preserved.
Here we content ourselves showing how the non-degeneracy conditions
\BHTa\ and \BHTb\ enter when solving the homological equation.

\medskip\noindent
At each iteration step we look for a transformation
$(\xi, \eta, \zeta, \sigma, \nu) \mapsto (x, y, z, \omega, \mu)$
with $\omega = \sigma+\Lambda_1(\sigma,\nu)$ and
$\mu = \nu+\Lambda_2(\sigma,\nu)$ independent from the variables
$(\xi, \eta, \zeta)$ so that the projection to the parameter space $P$ is preserved.
The transformation in the variables is generated by a $\sg$-equivariant
vector field $\Psi\in \cX^+$ that we write as
\[
   \Psi=U \partial_x + V \partial_{y} + W \partial_{z}.
\]
The homological equation reads
\begin{equation}
\label{homological}
   \ad N(X)(\Psi) = L+N
\end{equation}
with 
\[
   L_{\sigma,\nu}(\xi,\eta,\zeta)= \{Z-X\}_{lin,d} 
   \quad \mbox{and} \quad 
   N_{\sigma,\nu}(\xi,\eta,\zeta)=
   \Lambda_{1}(\sigma,\nu)\partial_{\xi} + 
   \Omega(\Lambda_2(\sigma,\nu))\zeta \partial_{\zeta}
\]
and determines the unknown $U$, $V$, $W$, $\Lambda_1$ and~$\Lambda_2$
according to
\begin{eqnarray}
\label{homological-sys}
   U_{\xi}\sigma & = &
   \Lambda_{1}+\widetilde{f}
\nonumber\\
   V_{\xi} \sigma + V_{\zeta} \Omega(\nu) \zeta & =& 
   \widetilde{g} +
   \widetilde{g}_{\eta} \eta + 
   \widetilde{g}_{\zeta} \zeta 
\\
   W_{\xi} \sigma + [ \Omega(\nu)\zeta,W] & = &
   \Omega (\Lambda_{2})\zeta + 
   \widetilde{h} + \widetilde{h}_{\eta} \eta +
   \widetilde{h}_{\zeta} \zeta . \nonumber
\end{eqnarray}
Here $U,V,W, \widetilde{f}, \widetilde{g}, \widetilde{h}$ and their derivatives
depend on $(\xi,0,0,\sigma,\nu)$. 
Moreover, greek subscripts denote derivatives,  
while $U_{\xi}\sigma = \Sigma_{j=1}^{n}U_{\xi_j}\sigma_j$ and similarly
for $V$ and~$W$.
These linear equations are solved by suitably truncated Fourier series.
Note that the left hand side of (\ref{homological-sys}) consists of the
components of the vector field $\ad N(X_{\sigma, \nu})(\Psi)$, where 
\[
   N(X_{\sigma,\nu})(\xi,\eta,\zeta)=
   \sigma\partial_{\xi}+\Omega(\nu)\zeta\partial_{\zeta}.
\]
For a given $Z$ (and hence $L$), the goal is to find $\Psi \in \cX^+_{lin,d}$
and $N\in \ker \ad N(X)^{T}\subseteq \cX^+_{lin,d}$ so that the homological
equation~\eqref{homological} is satisfied.
Here $\cX^+_{lin,d}=\cX^+_{lin}\cap \cX^+_{d}$ denotes the intersection set 
of the Taylor and Fourier truncations of vector fields in $\cX^+$. 

\medskip\noindent
We make the ansatz
\begin{equation}
   V(\xi,\eta,\zeta,\sigma,\nu)=V_0+V_1\eta+V_2\zeta 
   \quad \mbox{and} \quad 
   W(\xi,\eta,\zeta,\sigma,\nu)=W_0+W_1\eta+W_2\zeta 
\end{equation}
for the unknown $\Psi$, where $V_j$ and $W_j$ ($j=0,1,2$) depend on $\xi$
and on the multiparameter $(\sigma,\nu)$. 
Fourier expanding in $\xi$ and comparing coefficients in (\ref{homological-sys})
yields the following equations for an explicit (formal) construction of $\Psi$. 
To avoid clumsy notation we suppress the dependence on $(\sigma,\nu).$ 

\medskip\noindent
For $k\neq 0$, equation (\ref{homological-sys}) implies 
\begin{eqnarray}
   \I \langle k,\sigma \rangle U_k & = &  \widetilde{f}_k
\label{knonzerobegin}\\
   \I \langle k,\sigma \rangle V_{0,k} & = &\widetilde{g}_k,
\\
   \I \langle k,\sigma \rangle V_{1,k} & = & (\widetilde{g}_\eta)_k
\\
   V_{2,k} \left[\I \langle k,\sigma \rangle \id +
   \Omega(\nu)\right] & = &  (\widetilde{g}_\zeta)_k
\\
   \left[\I \langle k,\sigma \rangle \id -
   \Omega(\nu)\right] W_{0,k} & = & \widetilde{ h}_k
\\
   \left[\I \langle k,\sigma \rangle \id -
   \Omega(\nu)\right] W_{1,k} & = &  (\widetilde{h}_\eta)_k \label{knonzeroend}
\\
\label{MCC}
   \left[\I \langle k,\sigma \rangle \id - 
   \ad \Omega(\nu)\right] W_{2,k} & = & (\widetilde{h}_\zeta)_k 
\end{eqnarray}
and, similarly, for $k=0$  
\begin{eqnarray}
   - \Lambda_1& = & \widetilde{f}_0 \label{Lambda1}\\
   V_{2,0} \Omega(\nu) & = &  (\widetilde{g}_\zeta)_0 \label{V2} \\
   -\Omega(\nu) W_{0,0} & = &  \widetilde{h}_0 \label{W0}\\
   -\Omega(\nu) W_{1,0}& = &  (\widetilde{h}_\eta)_0 \label{W1}\\
\label{W2}
   - \ad \Omega(\nu) W_{2,k}  - \Omega(\Lambda_2) & = & (\widetilde{h}_\zeta)_0. 
\end{eqnarray}
On the one hand, it is clear by the Diophantine conditions
that for $k\neq 0$ none of the coefficients at the right hand
sides of~\eqref{knonzerobegin}-\eqref{knonzeroend}
is in the kernel, i.e.\ none of the eigenvalues 
$\I \langle k,\sigma \rangle$, 
$\I \langle k,\sigma \rangle\pm \lambda_j$, 
with $\lambda_j$ eigenvalue of $\Omega(\nu)$ are zero. 
For $k=0$, the equations (\ref{V2})-(\ref{W1}) are solvable by the
non-degeneracy condition~\BHTa\ since the right hand sides lie in~$\cB^-$.
The so-called solvability condition~\eqref{Lambda1} determines the
$\partial_{\xi}$-component
\begin{displaymath}
   \Lambda_{1}(\sigma, \nu)=-\frac{1}{(2\pi)^n}
   \int_{\T^n}\widetilde{f}(\xi,0,0,\sigma,\nu) \, {\rm d}\xi
\end{displaymath}
of~$N$ in~\eqref{homological}.
Turning our attention to equation (\ref{MCC}), we see that it admits the solution 
\[
   W_{2,k} = \left[ \I \langle k,
   \sigma \rangle \id -\ad \Omega(\nu) \right]^{-1}(\widetilde{h}_{\xi})_k
\]
if and only if the operator
$\left[\I \langle k,\sigma \rangle \id -\ad \Omega(\nu)\right]$
is invertible, which boils down to the condition 
\[
   \I \langle k,\sigma \rangle \neq \lambda_j-\lambda_l
\]
on the spectrum of $\ad\Omega(\nu)$, where $\lambda_j$
is an eigenvalue of~$\Omega(\nu)$. 
This inequality is the second Mel'nikov condition and
again guaranteed by the Diophantine conditions. 
For $k=0$ the splitting
\begin{equation}
\label{splitting}
   \im  (\ad_+ (\Omega_0 ))\oplus \ker (\ad_- (\Omega_0^T )) =  \gl_- (2p,\R),
\end{equation}
detailed in the Appendix lies at the basis of solving equation~\eqref{W2}. 
Indeed, the non-degeneracy condition \BHTb\ guarantees that we may choose
the \lcu\ for $\Omega$.
Using the Implicit Function Theorem and the fact that $\Omega$
(by construction) is an isomorphism 
between parameter spaces, it follows that (\ref{W2}) admits the solution 
\begin{equation}
   \Lambda_2(\sigma,\nu) =
   \Omega^{-1}\Bigl( - \pi \left(
   \widetilde{h}_{\zeta,0} + \ad \Omega(\nu)W_{2,0}
   \right) \Bigr), 
\end{equation}
where the mapping $\pi$ denotes the projection of $\gl_{-}(n,\R)$ onto the subspace
$ \ker (\ad_- (\Omega_0^T ))$ according to the splitting (\ref{splitting}).
Compare with \cite{CMC2}, Lemma~8.1.

\appendix

\section{Unfolding reversible linear matrices}
\label{appendix}

Let $\Omega_0\in \gl_-(2p;\R)$ be given; the aim of this appendix is to summarize
some results from \cite{ CMC2,Hov1, KnVdb1} which allow to describe a miniversal
unfolding of  $\Omega_0$, and to work out the details for two particular cases.

\medskip\noindent
Let $\Omega_0=\mathcal{S}_0 + \mathcal{N}_0$ be the Jordan-Chevalley 
decomposition of $\Omega_0$ into commuting semisimple and nilpotent parts.
The uniqueness of this decomposition implies that both $\mathcal{S}_0$ and 
$\mathcal{N}_0$ belong to $\gl_-(2p;\R)$.
Also
\begin{equation}
\label{ker-ad-Omega}
   \ker \ad(\Omega_0) = \ker \ad(\mathcal{S}_0) \cap \ker \ad(\mathcal{N}_0),
\end{equation}
as easily follows from the fact that $\mathcal{S}_0$ and $\mathcal{N}_0$ commute,
and as shown in~\cite{CMC2, KnVdb1} furthermore
\begin{equation}
\label{ker-ad-OmegaT}
   \ker \ad(\Omega_0^T) = \ker \ad(\mathcal{S}_0) \cap \ker \ad(\mathcal{N}_0^T).
\end{equation}
We know from (\ref{tangent-to-orbit}) that
$T_{\Omega_0}\orb(\Omega_0) = \im\,(\ad_+(\Omega_0))$,
while a classical result from linear algebra shows that
the subspace $\ker\,(\ad_-(\Omega_0^T))$ of $\gl_-(2p;\R)$
forms a complement of the tangent space $T_{\Omega_0}\orb(\Omega_0)$
to the orbit through $\Omega_0$.
Finally,
$\ker\,(\ad_-(\Omega_0^T))= \ker\,(\ad_-(\mathcal{S}_0)) \cap \ker\,(\ad_-(\mathcal{N}_0^T))$
by (\ref{ker-ad-OmegaT}), and hence we obtain the following result.

\begin{theorem}
\label{thm-lcu}
Let $\Omega_0\in\gl_-(2p;\R)$ be given, and let
$\Omega_0=\mathcal{S}_0 + \mathcal{N}_0$ be the Jordan-Chevalley
decomposition of $\Omega_0$. Then
\begin{equation}
\label{LCU-anyp}
   \Omega: \ker\,(\ad_-(\mathcal{S}_0)) \cap
   \ker\,(\ad_-(\mathcal{N}_0^T)) \longrightarrow \gl_-(2p;\R),
   \quad A\mapsto \Omega_0+A,
\end{equation}
forms a miniversal unfolding of $\Omega_0\in\gl_-(2p;\R)$.
\end{theorem}

\noindent
The unfolding $\Omega(\mu)$ is in the {centralizer} of $\mathcal{S}_0$. 
In the present context of linear systems one calls such an unfolding a
linear centralizer unfolding (\lcu\ for short).
Also note that $\Omega(\mu)-\Omega_0$ is linear in the unfolding parameters.
For the convenience of the reader we now explicitly work out a linear centralizer
unfolding~\eqref{LCU-anyp} for three particular choices of~$(\Omega_0, R)$.

\subsection{Unfolding multiple non-zero normal frequencies}

For our first example we assume that $\Omega_0\in\gl_-(2p;\R)$ has a
$1:1:\cdots:1$ resonance (or $p$-fold resonance), meaning that $\Omega_0$
has a pair of purely imaginary eigenvalues, say $\pm \I$, with algebraic
multiplicity $p$; we furthermore assume that we are in the generic situation,
with geometric multiplicity $1$.
The subspaces $\ker(\mathcal{N}_0^j)$ ($1\leq j\leq p$) form a strictly
increasing sequence of subspaces invariant under $\mathcal{S}_0$ and $R$, 
with $\dim \ker(\mathcal{N}_0^j) - \dim \ker(\mathcal{N}_0^{j-1})=2$. 
With respect to a conveniently chosen basis
$\{u_1^+,u_1^-,u_2^+,u_2^-,\ldots,u_p^+,u_p^-\}$
of $\R^{2p}$ the linear matrices $\Omega_0$ and $R$ have the matrix form
\begin{equation}
\label{eq:re-generic}
   \Omega_0 =\left(
   \begin{array}{ccccc}
      {\rm J}_2 &  {\rm J}_2 & {\rm O}_2 & \ldots & {\rm O}_2  \\
      {\rm O}_2 &  {\rm J}_2 & {\rm J}_2 & \ddots & \vdots\\
      \vdots & \ddots  & \ddots & \ddots & {\rm O}_2 \\
      \vdots & \   & \ddots & \ddots & {\rm J}_2\\
      {\rm O}_2 & \ldots  & \ldots & {\rm O}_2& {\rm J}_2 
   \end{array} 
   \right), \quad  
   R=\left(
   \begin{array}{ccccc}
      R_2 &  {\rm O}_2 & {\rm O}_2 & \ldots &{ \rm O}_2  \\
      {\rm O}_2 &  R_2 & {\rm O}_2 & \ddots & \vdots\\
      \vdots & \ddots  & \ddots & \ddots & \rm{O}_2 \\
      \vdots & \   & \ddots & \ddots & {\rm O}_2\\
      {\rm O}_2 & \ldots  & \ldots &{\rm O}_2& R_2 
   \end{array} 
   \right)
\end{equation}
with 
\begin{equation}
\label{J2O2R0}
   {\rm J}_2 = \left(%
   \begin{array}{cc}
      \phantom{-}0 & 1 \\
       -1 & 0 \\
   \end{array}
   \right), \quad
   {\rm O}_2 = \left(
   \begin{array}{cc}
      0 & 0 \\
      0 & 0 \\
   \end{array}
   \right) \quad \mbox{and} \quad
   R_2=\left(
   \begin{array}{cc}
      1 & \phantom{-}0\\
      0& -1
   \end{array}
   \right). 
\end{equation}
From this one can compute an \lcu\ of $\Omega_0$ as follows.

\begin{lemma}
Fix an $A\in \ker\,(\ad_-(\mathcal{S}_0)) \cap \ker\,(\ad_-(\mathcal{N}_0^T))$.
Then there exist constants $\mu_1,\mu_2,\ldots,\mu_p\in \R$
such that if we set
\[
A_j:= A - \sum_{i=1}^j\,\mu_i\mathcal{S}_0^i\left(\mathcal{N}_0^T\right)^{i-1},
\qquad (1\leq j \leq p),
\]
then $A_j(U_{p-i})=\{0\}$ for $1\leq j\leq p$ and $0\leq i\leq j-1$. 
In the particular case that $j=p$ we have 
\begin{equation}
\label{lcu-p-fold-res}
A= \sum_{i=1}^p \mu_i\mathcal{S}_0^i\left(\mathcal{N}_0^T\right)^{i-1}.
\end{equation}
\end{lemma}

\noindent
Combining (\ref{eq:re-generic}) and (\ref{lcu-p-fold-res}) an \lcu\ of $\Omega_0$
takes the explicit form
\begin{equation}\label{lcu-generic}
   \Omega(\mu) =\Omega_0+\left(
   \begin{array}{ccccc}
      \mu_1 {\rm J}_2 &  {\rm O}_2 & {\rm O}_2 & \ldots & {\rm O}_2  \\
      \mu_2 {\rm J}_2 &  \mu_1 {\rm J}_2 & {\rm O}_2 & \ddots & \vdots\\
      \vdots & \ddots  & \ddots & \ddots & {\rm O}_2 \\
      \vdots & \   & \ddots & \ddots & \rm{O}_2\\
      \mu_p\rm{J}_2 & \ldots  & \ldots & \mu_2\rm{J}_2& \mu_1\rm{J}_2 
   \end{array}
   \right)
\end{equation}
with unfolding parameters $\mu_1, \ldots, \mu_p\in \R$. 
This construction invariably leads to the same {\lcu},
we therefore speak from now on of {\it the} \lcu.
In case all eigenvalues of $\Omega_0\in\gl_-(2p;\R)$ are purely imaginary,
non-zero and with geometric multiplicity $1$, the \lcu\ of $\Omega_0$ can be
obtained by considering the different pairs of eigenvalues $\pm \I\alpha_j$,
multiplying (\ref{lcu-generic}) with $\alpha_j$ (using for each $j$ the
appropriate dimension and a new set of parameters), and juxtaposing the
obtained unfoldings as blocks along the diagonal.

\subsection{Unfolding multiple eigenvalue zero}

For our second and third example we assume that $\Omega_0\in\gl_-(2p;\R)$
has $0$ as an eigenvalue with geometric multiplicity $1$ and algebraic
multiplicity $2p$; then $\mathcal{S}_0=0$, $\mathcal{N}_0=\Omega_0$,
$\mathcal{N}_0^j\neq 0$ for $1\leq j <2p$, and $\mathcal{N}_0^{2p}=0$.
The subspaces $\ker(\mathcal{N}_0^j)$, $1\leq j \leq 2p$ are invariant under $R$;
they form a strictly increasing sequence,
with $\dim \ker(\mathcal{N}_0^j) - \dim \ker(\mathcal{N}_0^{j-1})=1$.
With respect to a conveniently chosen basis $\Omega_0=\mathcal{N}_0$ is a
classical nilpotent Jordan matrix with $1$'s above the diagonal.
The matrix form of $R$ depends on whether
$\ker(\mathcal{N}_0)\subset \textrm{Fix}(R)$, in which case $R$ has the same
matrix form as in (\ref{eq:re-generic}), or
$\ker(\mathcal{N}_0)\subset \textrm{Fix}(-R)$, whence the matrix form of $R$
equals minus the expression in (\ref{eq:re-generic}).

\medskip\noindent
To determine the \lcu\ of $\Omega_0$ we first consider some
$A \in \ker (\ad(\mathcal{N}_0^T))$; one easily shows that $A$ can be written as
$A= \sum_{j=1}^{2p}\nu_j\left(\mathcal{N}_0^T\right)^{j-1}$, with some constants $\nu_j\in\R$
($1\leq j\leq 2p$). Imposing the further condition that $A\in \gl_-(2p;\R)$ gives $\nu_j=0$ for
$j$ odd; setting $\mu_j:=\nu_{2j}$ for $1\leq j \leq p$ we obtain then the following \lcu\!\!:
\[
   \Omega(\mu) = \Omega_0 + \sum_{j=1}^p\,\mu_j\left(\mathcal{N}_0^T\right)^{2j-1}, \quad
   \mu=(\mu_1,\mu_2,\ldots,\mu_p)\in\R^p.
\]
Hence $\Omega_0$ has co-dimension $p$ and the \lcu\ is given by 
\begin{equation}
\label{unf_0_geomult_1}
   \Omega(\mu) = \left(
   \begin {array}{ccccccc}
      0&1&0&0&0&\cdots&0\\
      &0&1&0&0&\cdots&0\\
      & &0&1&0&\cdots&0\\
      & & &\ddots&\ddots&\ddots&\vdots\\
      & & & &\ddots&\ddots&0\\
      & & & & &\ddots&1\\
      & & & & & &0 
   \end{array}
   \right) + \left(
   \begin {array}{ccccccc}          
      0&\\
      \mu_1&0\\
      0&\mu_1&0\\
      \mu_2&0&\mu_1&0\\
      0&\mu_2&0&\ddots&\ddots\\
      \vdots&\ddots&\ddots&\ddots&\ddots&\ddots \\
      \mu_p &\ldots &0 &\mu_2& 0& \mu_1& 0
   \end{array}
   \right), 
\end{equation}
alternating diagonals with unfolding parameters $\mu_j$ and diagonals with~$0$.
Note that we may alternatively fix~$R$ to be of the form~\eqref{eq:re-generic}
and obtain the two cases by taking~\eqref{unf_0_geomult_1} and its transpose,
with $\Omega_0 = \mathcal{N}_0$ having its $1$'s below the diagonal.

\medskip\noindent
In case the condition $\dim \ker(\Omega_0)=1$ on the geometric multiplicity 
of the zero eigenvalue is dropped
the unfolding changes drastically and requires more parameters, i.e., has higher codimension. 
The same is true for our first example (non-zero normal frequencies). 
Further information on these cases can be found in \cite{Hov1}.


\begin{thebibliography}{99}
\bibitem{A71} V.I. Arnol'd. On matrices depending on parameters.
      {\it Russ.\ Math.\ Surv.}\ {\bf 26}(2), p. 29--43 (1971)
\bibitem{A83} V.I. Arnol'd.
      {\it Geometrical Methods in the Theory of
      Ordinary Differential Equations}.
      Springer (1983)
   \bibitem{bou97}
      J. Bourgain.
      On Melnikov's persistency problem.
      {\it Math.\ Res.\ Lett.}\ {\bf 4}, p. 445--458 (1997)
\bibitem{BB}
B.J.L. Braaksma and H.W. Broer.
      On a quasi-periodic Hopf bifurcation.
      {\it Ann.\ Inst.\ H.\ Poincar\'e, Analyse non lin\'eaire}
      {\bf 4}(2), p. 115--168 (1987)
\bibitem{BBH}
B.J.L. Braaksma, H.W. Broer and G.B. Huitema.
      Toward a quasi-periodic bifurcation theory.
      {\it Mem.\ AMS} {\bf 83} \#421, p. 83--175 (1990)
\bibitem{BrFu93}
      T.J. Bridges and J.E. Furter.
      {\it Singularity Theory and Equivariant Symplectic Maps}.
      LNM {\bf 1558}, Springer (1993)
\bibitem{BCH}
H.W. Broer, M.C. Ciocci and H. Han{\ss}mann.
      The quasi-periodic reversible Hopf bifurcation.
      In E. Doedel, B. Krauskopf and J. Sanders (eds.),
      {\it Recent Advances in Nonlinear Dynamics.
      Theme section dedicated to Andr\'e Vanderbauwhede},
      {\it Int.\ J. Bif.\ Chaos} {\bf 17}, p. 2605 -- 2623 (2007)
\bibitem{BGV2} H.W. Broer, M. Golubitsky and G. Vegter, 
Geometry of resonance tongues. 
In D. Ch\'eniot, N. Dutertre, C. Murolo, D. Trotman and A. Pichon (eds.), 
\textit{Proceedings of the 2005 Marseille Singularity School and Conference,
dedicated to Jean-Paul Brasselet on His 60th Birthday}, 
World Scientific, p. 327--356 (2007)
\bibitem{BHH07}
      H.W. Broer, H. Han{\ss}mann and J. Hoo.
      The quasi-periodic Hamiltonian Hopf bifurcation.
      {\it Nonlinearity} {\bf 20}, p. 417--460 (2007)
   \bibitem{BHJVW03}
      H.W. Broer, H. Han{\ss}mann, \`A. Jorba, J. Villanueva and
      F.O.O. Wagener.
      Normal-internal resonances in quasi-periodically forced
      oscillators: a conservative approach.
      {\it Nonlinearity} {\bf 16}, p. 1751--1791 (2003)
\bibitem{BHN} H.W.\ Broer, J.\ Hoo and V.\ Naudot.
Normal linear stability of quasi-periodic tori. 
     {\it J.\ Diff.\ Eq.}\ {\bf 232}(2), p. 355--418 (2007)
\bibitem{BH}
H.W. Broer and G.B. Huitema. Unfoldings of 
Quasi-periodic Tori in Reversible Systems. 
      {\it J.\ Dynamics Diff.\ Eq.}\ {\bf 7}, p. 191--212 (1995)
\bibitem{BHS} H.W.Broer, G.B. Huitema and M.B. Sevryuk. 
      {\it Quasi-Periodic Motions in Families of
      Dynamical Systems: Order amidst Chaos}.
      LNM {\bf 1645}, Springer (1996)
\bibitem{BHT}
H.W. Broer, G.B. Huitema and F.Takens. 
      Unfoldings of quasi-periodic tori.
      {\it Mem.\ AMS} {\bf 83} \#421, p. 1--81 (1990)
\bibitem{BS} 
H.W. Broer and M.B. Sevryuk.
     {\sc kam} Theory: quasi-periodicity in dynamical systems.
In: H.W. Broer, B. Hasselblatt and F. Takens (eds.) \textit{Handbook of Dynamical Systems,}
Volume 3. North-Holland (2008), to appear
   \bibitem{BrVe92}
      H.W. Broer and G. Vegter.
      Bifurcational Aspects of Parametric Resonance.
      {\it Dynamics Rep., new ser.}\ {\bf 1}, p. 1--53 (1992)
\bibitem{CMC2} 
M.C. Ciocci. {\em  Bifurcations of periodic orbits and persistence of quasi-periodic orbits 
in families of reversible systems.} PhD. Thesis, Universiteit Gent (2003)
\bibitem{CLB}
M.C. Ciocci, A. Litvak-Hinenzon and H.W. Broer, 
Survey on dissipative {\sc kam} theory including quasi-periodic bifurcation
theory based on lectures by Henk Broer.  
In:  J. Montaldi and T. Ratiu (eds.):
\textit{Geometric Mechanics and Symmetry: the Peyresq Lectures,} 
LMS Lecture Notes Series, {\bf 306}. 
Cambridge University Press, p. 303--355 (2005)
\bibitem{Gib} C.G. Gibson. {\it Singular points of smooth mappings}.  
Research Notes in Mathematics {\bf 25}, Pitman (1979)
\bibitem{h1h98}
 H. Han{\ss}mann.
 The Quasi-Periodic Centre-Saddle Bifurcation.
 \textit{J. Diff.\ Eq.}\ \textbf{142}, p. 305--370 (1998)
   \bibitem{h1h07}
      H. Han{\ss}mann.
      {\it Local and Semi-Local Bifurcations in Hamiltonian Dynamical Systems
      -- Results and Examples}.
      LNM {\bf 1893}, Springer (2007)
\bibitem{Hoo} J. Hoo.
      {\it Quasi-periodic bifurcations in a strong resonance: combination tones in gyroscopic stabilization}.
      PhD thesis, Rijksuniversiteit Groningen (2005)
\bibitem{Hov1}
I. Hoveijn.
      Versal Deformations and Normal Forms for Reversible and
      Hamiltonian Linear Systems.
      {\it J. Diff.\ Eq.}\ {\bf 126}(2), p. 408--442 (1996)
\bibitem{H} G.B. Huitema. {\it Unfoldings of quasi-periodic tori}. PhD thesis, Rijksuniversiteit Groningen (1988)
\bibitem{Iooss} G.\ Iooss.
A codimension $2$ bifurcation for reversible vector fields.
In W.F.~Langford and W.~Nagata (eds.):
      {\it Normal Forms and Homoclinic Chaos, Waterloo 1992}.
      {\ Fields Institute Communications} {\bf 4}, AMS, p. 201--217 (1995)
\bibitem{KnVdb} 
J. Knobloch and A. Vanderbauwhede.
A General Reduction Method for Periodic Solutions in Conservative and
Reversible Systems.
      {\it J.\ Dynamics Diff.\ Eq.}\ {\bf 8}, p. 71--102 (1996)
\bibitem{KnVdb1} 
J. Knobloch and A. Vanderbauwhede.
Hopf bifurcation at k-fold resonances in reversible systems.
Preprint Technische Universit\"at Ilmenau, No.\ M 16/95 (1995) \\
URL: \verb+http://cage.ugent.be/~avdb/articles/k-fold-Res-Rev.pdf+
\bibitem{lam94}
   J.S.W. Lamb.
   \emph{Reversing symmetries in dynamical systems}.
   PhD thesis, Universiteit van Amsterdam (1994)
\bibitem{LaCa95}
   J.S.W. Lamb and H.W. Capel.
Local bifurcations on the plane with reversing point group symmetry.
{\it Chaos, Solitons and Fractals} {\bf 5}(2), p. 271--293 (1995)
   \bibitem{mel65}
      V.K. Mel'nikov.
      On some cases of conservation of conditionally periodic motions
      under a small change of the Hamiltonian function.
      {\it Sov.\ Math.\ Dokl.}\ {\bf 6}(6), p. 1592--1596 (1965)
\bibitem{Mos 67} J.K. Moser.
      Convergent Series Expansion for Quasi-Periodic Motions.
      {\it Math.\ Ann.}\ {\bf 169}(1), p. 136--176 (1967)
\bibitem{Moser} J.K. Moser.
{\it Stable and Random Motions in Dynamical Systems.}
Ann.\ Math.\ Stud.\ {\bf 77}, Princeton Univ.\ Press (1973)
\bibitem{Po82}J. P\"oschel.
      Integrability of Hamiltonian Systems on Cantor Sets.
      {\it Comm.\ Pure Appl.\ Math.}\ {\bf 35}, p. 653--696 (1982)
\bibitem{PCQW90}
T. Post, H.W. Capel, G.R.W. Quispel, and J.P. van der Weele.
Bifurcations in two-dimensional reversible maps.
{\it Physica A} {\bf 164}(3), p. 625--662 (1990)
\bibitem{rue01}
H. R\"ussmann.
Invariant tori in non-degenerate nearly integrable Hamiltonian systems.
{\it Reg.\ \& Chaot.\ Dyn.}\ {\bf 6}(2), p. 119--204 (2001)
\bibitem{sev86}
   M.B.\ Sevryuk.
   \emph{Reversible systems}.
   LNM {\bf 1211}, Springer (1986)
\bibitem{Sev1} 
M.B. Sevryuk. Linear reversible systems and their versal deformations.
{\it J. Soviet Math.}\ {\bf 60}(5), p. 1663--1680 (1992) 
\bibitem{sev95} 
M.B. Sevryuk.
The iteration-approximation decoupling in the reversible KAM~theory.
{\it Chaos} {\bf 5}(3), p. 552--565 (1995)
\bibitem{wag05}
      F.O.O. Wagener.
      On the quasi-periodic $d$-fold degenerate bifurcation.
      {\it J. Diff.\ Eq.}\ {\bf 216}, p. 261--281 (2005)
\bibitem{Ba}
   Wei Baoshe.
   Perturbations of lower dimensional tori in the resonant zone
   for reversible systems.
   \emph{JMAA} {\bf 253}(2), p. 558--577 (2001)
   \bibitem{XuYo01}
      J. Xu and J. You.
      Persistence of lower-dimensional tori under the first Melnikov's
      non-resonance condition.
      {\it J. Math.\ Pures Appl.}\ {\bf 80}(10), p. 1045--1067 (2001)
\end{thebibliography}
\end{document}